\documentclass{article}

\usepackage{arxiv}

\usepackage{amsmath, amsthm, amssymb}
\usepackage{booktabs}
\usepackage{dsfont}
\usepackage{algorithm}
\usepackage{algpseudocode}
\usepackage{graphicx}
\usepackage{hyperref}
\usepackage{mathtools}
\usepackage{multirow}
\usepackage[square,numbers]{natbib}
\usepackage{pifont}
\usepackage{subcaption}
\usepackage[table, svgnames, dvipsnames]{xcolor}
\usepackage{tikz}
\usepackage{url}

\usepackage{algorithm}
\usepackage{algpseudocode}

\definecolor{prettyblue}{rgb}{0.2, 0.2, 0.6}
\definecolor{cocoabrown}{rgb}{0.82, 0.41, 0.12}
\definecolor{yesgreen}{rgb}{0, 153, 51}

\newcommand{\vb}{\boldsymbol b}
\newcommand{\vc}{\boldsymbol c}

\newcommand{\vv}{\boldsymbol v}

\newcommand{\vx}{\boldsymbol x}

\newcommand{\mA}{\boldsymbol A}

\newcommand{\mC}{\boldsymbol C}

\newcommand{\mV}{\boldsymbol V}
\newcommand{\mW}{\boldsymbol W}

\newcommand{\trans}{^\mathsf{T}}

\title{Learning Optimal Objective Values for MILP}
\date{}

\author{
 Lara Scavuzzo \\
  TU Delft\\
  \texttt{l.v.scavuzzomontana@tudelft.nl} \\
   \And
 Karen Aardal \\
  TU Delft\\
  \texttt{k.i.aardal@tudelft.nl} \\
  \And
 Neil Yorke-Smith \\
    TU Delft\\
    \texttt{n.yorke-smith@tudelft.nl} \\
}

\begin{document}
\maketitle
\begin{abstract}
Modern Mixed Integer Linear Programming (MILP) solvers use the Branch-and-Bound algorithm together with a plethora of auxiliary components that speed up the search. In recent years, there has been an explosive development in the use of machine learning for enhancing and supporting these algorithmic components \cite{Scavuzzo2024}. Within this line, we propose a methodology for predicting the optimal objective value, or, equivalently, predicting if the current incumbent is optimal. For this task, we introduce a predictor based on a graph neural network (GNN) architecture, together with a set of dynamic features. Experimental results on diverse benchmarks demonstrate the efficacy of our approach, achieving high accuracy in the prediction task and outperforming existing methods. These findings suggest new opportunities for integrating ML-driven predictions into MILP solvers, enabling smarter decision-making and improved performance.
\end{abstract}



\section{Introduction}
\label{sec:intro}

Mixed Integer Linear Programming (MILPs) is a widespread tool for modelling mathematical optimization problems, with applications in numerous real-world scenarios. The Branch-and-Bound (B\&B) algorithm, which employs a divide-and-conquer approach, is the preferred method for solving MILPs to global optimality.  In recent years, there has been a surge in interest in harnessing the power of machine learning (ML) tools to aid the solution process of MILPs.  From solution prediction (e.g. \cite{Ding2020, Nair2020, Sonnerat2021}) to interventions on the heuristic rules used by the solvers (e.g. \cite{Gasse2019, Chmiela2021, Paulus2022}), several approaches have been studied in the literature (see \citet{Scavuzzo2024} for an in-depth discussion of this topic). The overarching trend is to build dynamic MILP solvers that can make active use of the large amounts of data produced during the solving process. 

Many of the decisions that must be made during the B\&B process could be better informed were the optimal solution known from the start. In fact, even knowing the optimal \emph{objective value} can positively influence the solver behaviour. For example, once a solution is found that matches this value, any effort to find new solutions can be avoided. With perfect information of the optimal objective value, a solver can further do more aggressive pruning of nodes. In general, having this knowledge can allow the solver to adapt its configuration, putting more emphasis on different components. Even in absence of perfect information, a good prediction of the optimal objective value can still be used to change the solver settings or to devise smarter rules, such as node selection policies that account for this predicted value. Inspired by these observations we ask the two following closely-related questions:
\begin{itemize}
    \item[\textbf{(Q1)}] How well can we predict the optimal \emph{objective value}?
    \item[\textbf{(Q2)}] With what accuracy can we predict, during the solution process, whether or not a given solution is optimal?
\end{itemize}

Our contributions are as follows. First, we propose a methodology to predict optimal objective values, answering (Q1). We then use the output of this predictive model, together with additional data, as input of our proposed classifiers, which give an answer to question (Q2). For this second task, we also propose some metrics that capture the state of the solving process, and that prove to be valuable for our classifier.
Our computational study shows the high accuracy of our proposed predictor. Furthermore, when compared to previous methods, our classifiers show better performance. Finally, we provide further insight into how the performance can be tuned to the desired behaviour and into the ways that the classifier makes use of the provided data.

Our discussion is organized as follows. We start by defining some key concepts and notation in Section~\ref{sec:background}, followed by a discussion of the work most closely related to ours (Section~\ref{sec:related}). Section~\ref{sec:methodology} describes our methodology in detail. The results of our computational study are presented in Section~\ref{sec:results}. Finally, we conclude with some final remarks and future work in Section~\ref{sec:conclusions}. The code to reproduce all experiments is available online \cite{Git-objvalprediction}.

\section{Background}
\label{sec:background}
\paragraph{Mixed Integer Linear Programming}
Given are a matrix $\mA\in \mathbb{Q}^{m\times n}$, vectors $\vc\in\mathbb{Q}^n$ and $\vb\in\mathbb{Q}^m$, and a partition $(\mathcal{I}, \mathcal{C})$ of the variable index set $\{1,2,...,n\}$. A Mixed Integer Linear Program is the problem of finding
\begin{equation}
\label{eq:MILP}
    \begin{aligned}
 z^* = \min\ \  & \vc^{\trans} \vx  \\
 \text{subject to } & A \vx \geq \vb ,  \\
 & x_{j}\in \mathbb{Z}_{\geq 0}\quad \forall j\in \mathcal{
 I}, \\
 & x_{j}\geq 0 \quad \forall j\in \mathcal{C}. \\
\end{aligned}
\end{equation}
Notice that the variables in $\mathcal{I}$ are required to be integer. Removing this integrality constraint turns the problem into a Linear Program (LP), which constitutes a relaxation of the original MILP, known as the \emph{LP relaxation}. While MILP is $\mathcal{NP}$-hard, LPs are polynomial solvable.

\paragraph{Solving Mixed Integer Linear Programs}
The standard approach to solving MILPs is to use the LP-based branch-and-bound (B\&B) algorithm. This algorithm sequentially partitions the feasible region, while using LP relaxations to obtain lower bounds on the quality of the solutions of each sub-region. This search can be represented as a binary\footnote{Standard implementations of the B\&B algorithm use single-variable disjunctions that partition the feasible set into two. Other approaches exist but are, to the best of our knowledge, not implemented in standard optimization software.} tree. At a given time $t$ of the solution process we use $T_t$ to denote the search tree, i.e. the set of nodes, constructed so far by the B\&B algorithm. We denote by $\vx^*$ the optimal solution to Problem (\ref{eq:MILP}) and $z^*$ its corresponding optimal objective value. For a given node $i$ of the search tree, let $z^{LP}_i$ be the optimal objective value of the node's LP relaxation. We use the notation $z^{LP}$ for the root node, i.e., the solution to the original problem's LP relaxation. At any point of the search, an integer feasible solution provides an upper bound on the optimal objective value. Let $\bar{\vx}(t)$ be the best known solution at time $t$ and let $\bar{z}(t)=\vc^{\trans} \bar{\vx}(t)$ denote its objective value (also called the \emph{incumbent}). Then we can prune any node $i$ such that $z^{LP}_i \geq \bar{z}(t)$. 

The nodes of $T_t$ can be classified into three types:
\begin{itemize}
    \item $I_t$ is the set of inner nodes of the tree. This is, nodes that have been processed (its LP relaxation solved) and resulted in branching.
    \item $L_t$ is the set of leaves of the tree. This is, the set of nodes that have been processed and resulted in pruning or in an integer feasible solution.
    \item $O_t$ is the set of open nodes, i.e., nodes that have not been processed yet.
\end{itemize}

As mentioned before, the incumbent $\bar{z}(t)$ provides an upper bound on $z^*$. We can also obtain a global lower bound.  Let $\underline{z}(t):= \min_{i\in O_t}\{z^{LP}_i\}$. Then notice that necessarily $\underline{z}(t) \leq z^*$.

In practice, MILP solvers implement a plethora of other techniques to accelerate the solution process. Among them, cutting planes and primal heuristics are essential parts of today's mathematical optimization software. 

\paragraph{MILP solving phases}
The B\&B algorithm can solve MILPs to optimality. This means that, if the algorithm terminates, it does so after having obtained a feasible solution and a proof of its optimality (or, on the contrary, proof of infeasibility). Several solver components work together for this goal, each with more or less focus on the feasibility and the optimality parts. \citet{Berthold2018} point out that, typically, the optimal solution is found well before the solver can prove optimality. Following this, they propose partitioning the search process into phases, according to three target goals. These phases are the following.
\begin{enumerate}
    \item \textbf{Feasibility.} This phase encompasses the time spanned from the beginning of the search until the first feasible solution is found.
    \item \textbf{Improvement.} From the moment the first feasible solution is found until an optimal solution is found.
    \item \textbf{Proving.} Spans the time elapsed from the moment the optimal solution is found until the solver terminates with a proof of optimality.
\end{enumerate}

The transition between the first and the second phase happens when a feasible solution is found. In contrast, the moment in which the solver transitions from the second to the third phase is unknown until the search is completed. Notice that if the instance is infeasible the solver terminates while in the first phase.  For the purpose of our study we assume that the instances are feasible.

\section{Related Work}
\label{sec:related}
\paragraph{MILP solution prediction}
In recent years, the topic of solution prediction for combinatorial optimization problems has gained momentum \citep{Shen2022AdaptiveSP,DBLP:conf/cpaior/EfthymiouY23,Han2023AGP}.  For MILPs, the goal is to produce a (partial) assignment of the integer variables via a predictive machine learning model. This prediction can then be used to guide the search in different ways. \citet{Ding2020} impose a constraint that forces the search to remain in a neighbourhood of the predicted optimal solution. In this way, by restricting the size of the feasible region, the authors aim to accelerate the solution process. In contrast, the approaches of \citet{Nair2020} and \citet{Khalil2022mipgnn} consist in fixing a subset of variables to their predicted optimal value, letting the solver optimize over the remaining ones. \citet{Khalil2022mipgnn} further propose a solver mode that uses the predicted solution to guide the node processing order. In the present work, we take a different path by aiming to predict the \emph{optimal objective value}, as opposed to the solution, i.e., the values that each variable takes. This task is easier from a learning perspective, yet still offers several ways in which one can exploit this information.

\paragraph{Phase transition predictions}
 \citet{Berthold2018} defined the three phases of MILP solving that were introduced in Section \ref{sec:background}. Their goal is to adapt the solver's strategy depending on the phase. For this purpose, they propose two criteria that can be used to predict the transition between phase 2 (improvement) and phase 3 (proving) without knowledge of the optimal solution. 
These criteria are based on node estimates: for every node $i\in T_t$, the solver SCIP keeps an estimate $\hat{c}(i)$ of the objective value of the best solution attainable at that node (see \cite{Berthold2018} for a formal definition of how this is computed). At time $t$ of the solving process, let $\hat{c}^{\min}(t):=\min \{\hat{c}(i) \mid i\in O_t\}$ be the minimum of these estimates among the open nodes. We further define $d(i)$ to be the depth\footnote{We define the depth of a node as its distance to the root node. Therefore, by definition, the depth of the root node is zero.} of node $i$. The first transition criterion, the \emph{best-estimate} criterion, indicates that the transition moment is the first time the incumbent becomes smaller than $\hat{c}^{\min}(t)$. Formally, let us define a binary classifier $C^{\text{est}}$ that indicates if the transition has occurred using the criterion

 \begin{equation} \label{eq:estimate_classifier}
C^{\text{est}} = 
\begin{cases}
        1 & \text{if } \min_{s\in [0,t]} \{\bar{z}(s) -
        \hat{c}^{\min}(s)\} < 0 \\
        0 & \text{otherwise}. 
\end{cases}
\end{equation}

The second criterion is called \emph{rank-1} and is based on the set of open nodes with better estimate than the processed nodes at the same depth. Formally, let
 $$R^1(t):=\Big\{i \in O_t \mid \hat{c}(i) \leq \inf \{\hat{c}(j):j\in I_t\cup L_t, d(j)=d(i)\} \Big\} \,$$

This set can be used to define a classifier $C^{\text{rank-1}}$ that indicates that the transition has occurred once the set becomes empty for the first time. This is,

 \begin{equation} \label{eq:rank1_classifier}
C^{\text{rank-1}} = 
\begin{cases}
        1 & \text{if } \min_{s\in [0,t]} |R^1(s)| = 0 \\
        0 & \text{otherwise}. 
\end{cases}
\end{equation}

The authors use these criteria to switch between different pre-determined solver settings depending on the phase of solving. Their experiments show improved solving time, especially when using the rank-1 criterion. However, it is also clear that both criteria tend to be satisfied \emph{before} the phase transition actually occurs, and there is some room for improvement in the accuracy of the classifiers, as we shall see from our own computational study.

\paragraph{B\&B resolution predictions}
Closely related to the present work is that of \citet{Hendel2022b}, who use a number of solver metrics to predict the final B\&B tree size. They use a combination of metrics from the literature, together with their own, as input to a machine learning model that estimates the final tree size dynamically as the tree is being constructed. Their method was incorporated into version 7.0 of the solver SCIP as a progress metric for the user.
In a similar fashion, \citet{Fischetti2019} use a number of solver metrics to predict, during the solving process, whether or not the run will end within the given time limit. This prediction can be used to adapt the solver behaviour in the case that the answer is negative.

\section{Methodology}
\label{sec:methodology}

This section details the methodology used to answer questions Q1 (Section~\ref{sec:methodology_objpred}) and Q2 (Section~\ref{sec:methodology_phase}). We assume we are given a space $\mathcal{X}$  of instances of interest. For some tasks, we will use the bipartite graph representation of MILPs introduced by \citet{Gasse2019}. This is, given an MILP instance $X\in \mathcal{X}$ defined as in Eq.~\ref{eq:MILP}, we build a graph representation as follows: each constraint and each variable have a corresponding representative node. A constraint node is connected to a variable node if the corresponding variable has a non-zero coefficient in the corresponding constraint. Each node has an associated vector of features that describes it. We utilize the same features as Gasse et al., except that we do not include any incumbent information. In short, instead of the raw data in $X\in\mathcal{X}$ we use the graph representation, which we denote $X_G\in\mathcal{X}_G$, and is composed of a tuple $X_G=(\mC, \mV, \mA)$, where $\mC\in \mathbb{R}^{m\times d_c}$ and $\mV\in \mathbb{R}^{n\times d_v}$ represent the constraint and variable features, respectively, and  $\mA\in \mathbb{R}^{m\times n}$ is the adjacency matrix. \\

\begin{figure}
    \centering
    \includegraphics[width=\textwidth]{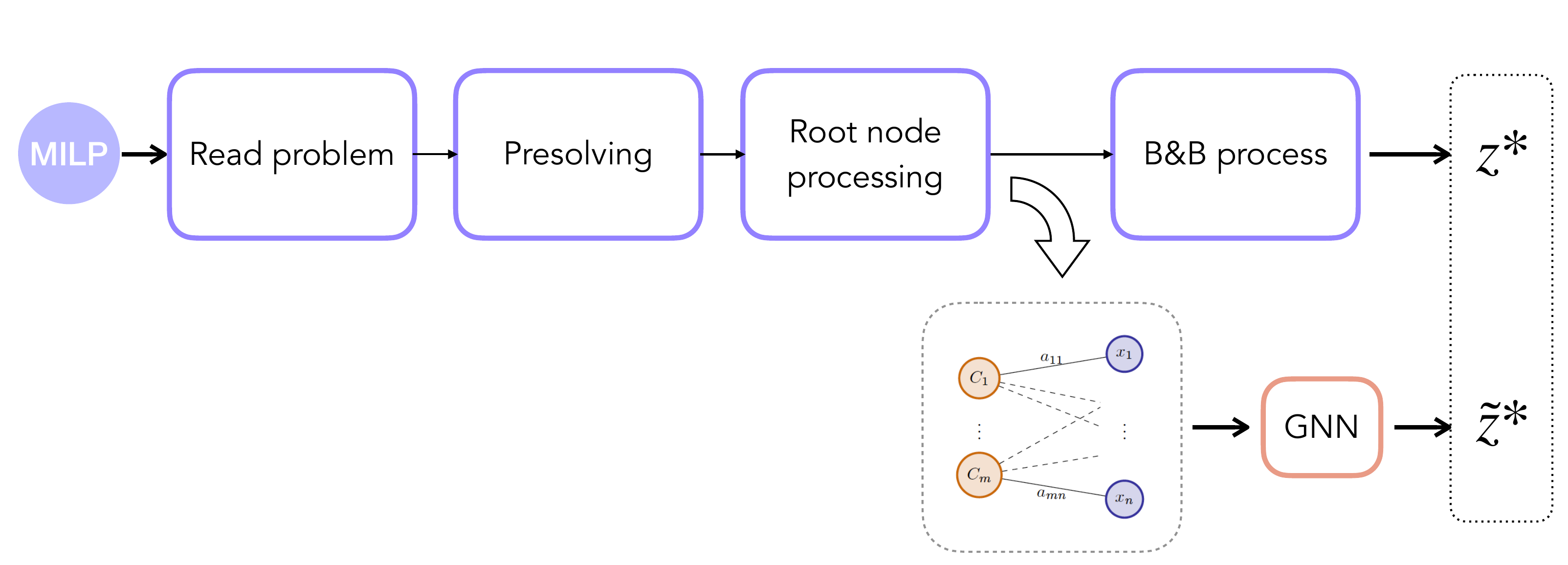}
    \caption{Optimal objective value prediction task. The MILP representation is computed after the root node has been processed. This serves as an input to a GNN that outputs a prediction $\tilde{z}^*$ of the optimal objective value.}
    \label{fig:process}
\end{figure}

\subsection{Optimal value prediction}
\label{sec:methodology_objpred}
The first task we tackle is the one of predicting the optimal objective value (Q1). That is, given an MILP instance $X \in \mathcal{X}$, we want to predict the optimal objective value $z^*$. This prediction is computed  once and for all at the root node, once the LP solution is available. We frame this as a regression task. This process is depicted in Figure~\ref{fig:process}. \\

For this regression task, we utilise the bipartite graph representation of \citet{Gasse2019} defined above, which is processed using a Graph Neural Network (GNN) that performs two half-convolutions. In particular, the feature matrices $\mC$ and $\mV$ first go through an embedding layer with two feedforward networks with ReLU activation. Next, one first pass updates the constraint descriptors using the variable descriptors, while a second pass updates the variable descriptors using the (new) constraint descriptors. This is done with message-passing operations, computed as
\begin{equation}
    \vc'_i = \mW^{(11)}\vc_i + \mW^{(12)}\sum_{j=1}^n \mA_{ij} \vv_j
\end{equation}
\begin{equation}
    \vv'_j = \mW^{(21)}\vv_j + \mW^{(22)}\sum_{i=1}^m \mA_{ij} \vc'_i
\end{equation}
where $\mW^{(11)}$, $\mW^{(12)}$, $\mW^{(21)}$ and $\mW^{(22)}$ are trainable weights, $\vc_i$ is the feature vector of constraint $i$ and $\vv_j$ is the feature vector of variable $j$. The variable descriptors then go through another feedforward network with ReLU activation. Finally, average pooling is applied to obtain one single output value.  

Our goal is to learn a mapping $f(X_g): \mathcal{X}_G \mapsto \mathbb{R}$ which outputs an approximation $\tilde{z}^*$ of the optimal objective value $z^*$. At the moment of this prediction, the solution to the root LP relaxation is known and can be used for further context. In order to exploit that knowledge, we test three potential targets for the machine learning model, namely 
$$\Theta_1 = z^*$$
$$\Theta_2 = \frac{z^*}{z^{LP}}$$
$$\Theta_3 = z^* -  z^{LP} \, .$$
This gives rise to three models $f_1(X_g)$, $f_2(X_g)$ and $f_3(X_g)$, which we later transform into the desired output by setting either $f(X_g)=f_1(X_g)$, $f(X_g)=f_2(X_g)\cdot z^{LP}$, or $f(X_g)=f_3(X_g) + z^{LP}$.

\subsection{Prediction of phase transition}
\label{sec:methodology_phase}
The second task (Q2) is predicting the transition between phases 2 (improvement) and 3 (proving). That is, at any point during the solution process we want to predict whether the incumbent is in fact optimal. We cast this problem as a classification task. \\

We test the performance of two classifiers. The first one is based on the output of the GNN model discussed in Section \ref{sec:methodology_objpred}. Given an instance $X\in\mathcal{X}$ (in fact its associated graph representation $X_G$) and the current incumbent $\Bar{z}$, we obtain a binary prediction $C^{GNN}_{\epsilon}:\mathcal{X}_G\times \mathbb{R} \mapsto \{0,1\}$ in the following way
\begin{equation} \label{eq:gnn_classifier}
C^{GNN}_{\epsilon}(X_G, \Bar{z}) = 
\begin{cases}
        1 & \text{if } \Bar{z} < f(X_G) + \epsilon \cdot |f(X_G)| \\
        0 & \text{otherwise} 
\end{cases}
\end{equation}
for some $\epsilon\in [-1,1]$. The $\epsilon$ parameter allows us to control the confidence in the prediction.

The $C^{GNN}_{\epsilon}$ classifier is static, in the sense that it does not make use of any information coming from the B\&B process. On the contrary, the second predictor we propose, which we call $C^D$, is based on a set of dynamic metrics that are collected during the solving process. The metrics are the following.

\paragraph{Gap} Following SCIP \citep{SCIP8}, we define the gap as
\begin{equation}
\label{eq:gap}
    g(t) :=
    \begin{cases*}
      1 & if no solution has been 
      found yet or $\bar{z}(t)\cdot \underline{z}(t)<0$,\\
      \frac{|\bar{z}(t)-\underline{z}(t)|}{\max \{|\bar{z}(t)|, |\underline{z}(t)|, \epsilon\}} & otherwise.
    \end{cases*}
\end{equation}

\paragraph{Tree weight} For a given node $v\in T_t$, let $d(v)$ denote the node's depth. Then, the tree weight at time $t$ is defined as
\begin{equation}
\label{eq:tree_weight}
    \omega(t) := \sum_{v\in L_t} 2^{-d(v)}  \, .
\end{equation}
This metric was first defined by \citet{Kilby2006}.


\paragraph{Median gap } Let $m(t) = \text{median}\{z^{LP}_i \mid i\in O_t\}$ and let $\bar{z}^0$ be the first incumbent found. We define the median gap as
\begin{equation}\label{eq:median_gap}
    \mu (t) = \frac{|\bar{z}(t) - m(t)|}{|\bar{z}^0- z^{LP}|}
\end{equation}

\paragraph{Trend of open nodes} For a certain window size $h$, we store the values of $|O_k|$ for $k\in \{t-h, t-h+1, ..., t\}$. We then fit a linear function using least squares to compute the \emph{trend} of this sequence. We denote this trend at time $t$ as $\tau(t)$.

\paragraph{Ratio to GNN prediction} We make use of the prediction $f(X_G)$ coming from the GNN model and include the ratio with respect to the current incumbent as a metric. In particular we use
\begin{equation}
    \rho (t) = \frac{f(X_G)}{\bar{z}(t)}
\end{equation}

Notice that, while the gap and the tree weight are metrics from the literature, the other three are our own.

The input to the classifier is therefore a tuple $X_D=(g(t), \omega(t), \mu(t), \tau(t), \rho(t) )$. We train a classifier $C^D(X_D)$ that makes use of these dynamic features to make a binary prediction on whether we are in phase 2 or 3. We use a simple logistic regression, which will allow us to more easily interpret the resulting model, in contrast to more complex machine learning models.

\section{Computational Results}
\label{sec:results}
This section describes our computational setup and results. All experiments were performed with the solver SCIP v.8.0 \cite{SCIP8}. Code for reproducing all experiments in this section is available online \cite{Git-objvalprediction}.

\subsection{Experimental Set Up}

\paragraph{Benchmarks}
We use three NP-hard problem benchmarks from the literature: set covering, combinatorial auctions and generalized independent set problem (GISP). We create a fourth benchmark (mixed) that is comprised of instances of the three types, in equal proportion. The method and configuration used for generation of the instances is summarized in Table \ref{tab:benchmarks}. For each instance type, we generate 10,000 instances for training, 2000 instances for validation and another 2000 for testing.

\begin{table*}[b]
\caption{Method and configuration settings used to generate the instances of problem benchmark.}
\label{tab:benchmarks}
\centering
\begin{center}
\footnotesize
\begin{tabular}{ccc}
    Benchmark & Generation method & Configuration \\
    \toprule
    \multirow{2}{*}{Set covering} & \multirow{2}{*}{\citet{Balas1980}} & Items: 750 \\ 
    & & Sets: 1000 \\
    \midrule
    Combinatorial & \citet{Leyton2000} & Items: 200 \\
    auctions & with arbitrary relationships & Bids: 1000 \\
    \midrule
    \multirow{4}{*}{GISP} &  & Nodes: 80 \\
     & \citet{Colombi2017} & $p=0.6$ \\
     & with Erdos-Renyi graphs & $\alpha=0.75$ \\ 
      & & SET2, A \\
    \bottomrule
\end{tabular}

\end{center}
\end{table*}

\begin{figure}
     \centering
     \begin{subfigure}[b]{0.45\textwidth}
         \centering
         \includegraphics[width=\textwidth]{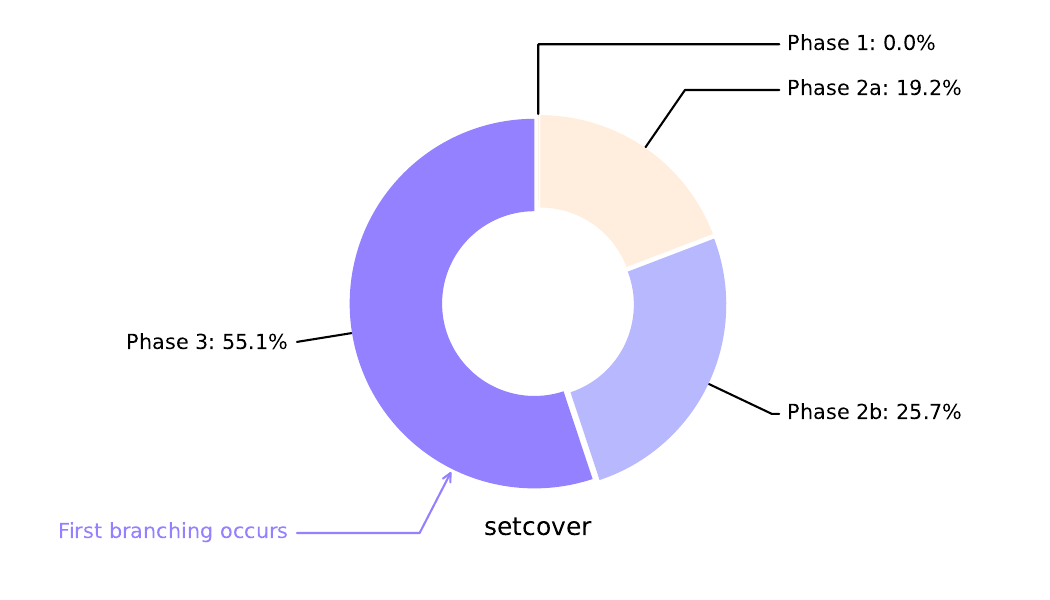}
         \caption{Set covering}
     \end{subfigure}
     \begin{subfigure}[b]{0.45\textwidth}
         \centering
         \includegraphics[width=\textwidth]{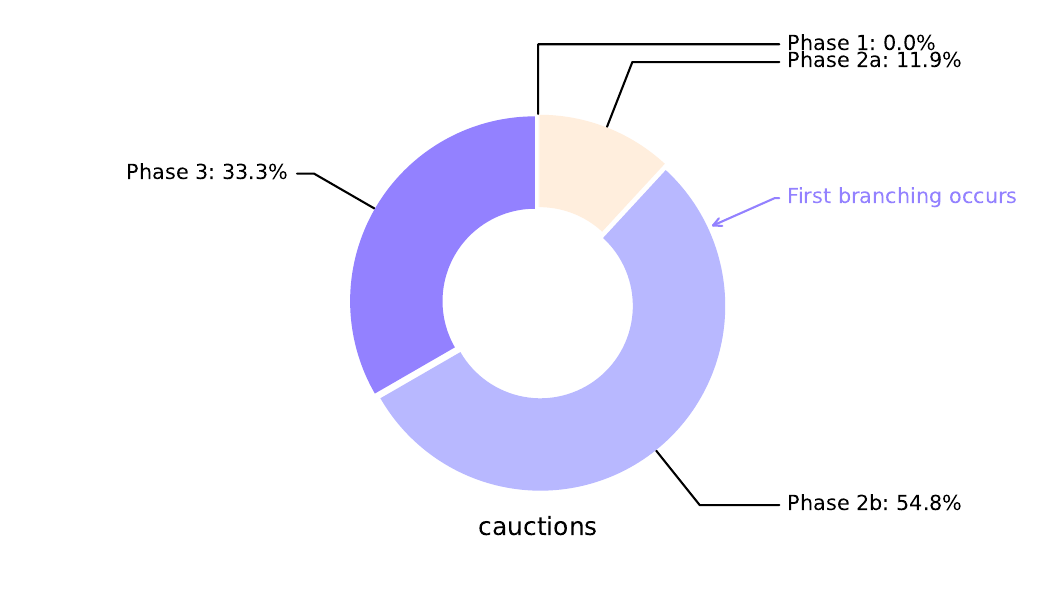}
         \caption{Combinatorial auctions}
     \end{subfigure} \\
     \begin{subfigure}[b]{0.45\textwidth}
         \centering
         \includegraphics[width=\textwidth]{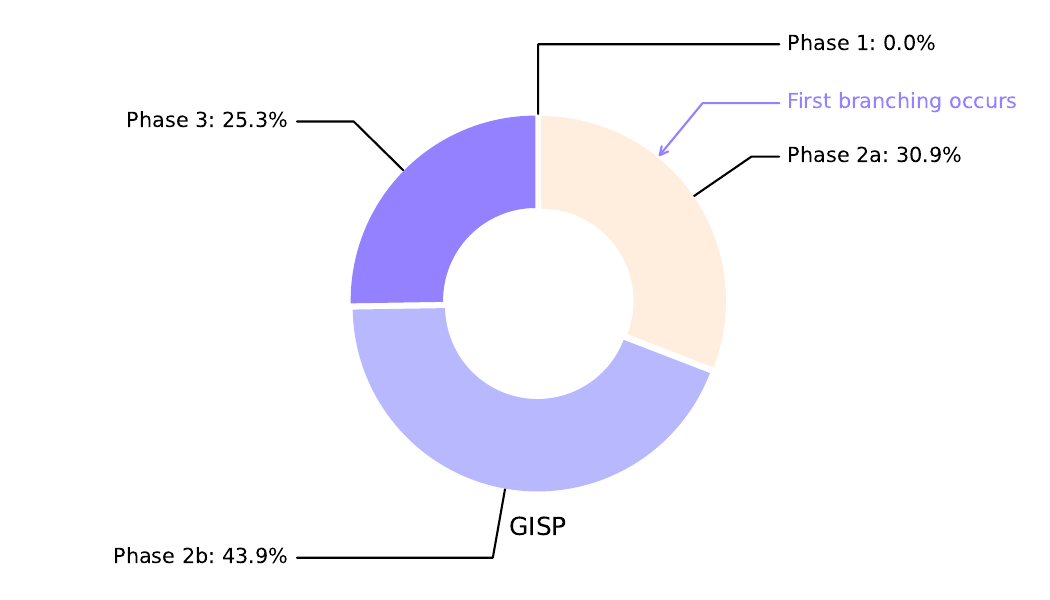}
         \caption{GISP}
     \end{subfigure}
\caption{Phase analysis of three instance types. We divide the solution process into (1) Feasibility, in dark yellow, (2a) Improvement up to 5\% to optimality, in light yellow, (2b) Improvement from 5\% to optimal, in light purple, and (3) Proving, in dark purple. We also indicate when the first branching occurs. The data is averaged over 100 instances with 3 randomization seeds (i.e., 300 samples).}
\label{fig:phase_analysis}
\end{figure}

\paragraph{Phase analysis}
As a first approach to the instances, we run an experiment to analyze the breakdown into solving phases. We solve 100 of the training instances, each with 3 different randomization seeds, which gives us a total of 300 data points per benchmark. During the solution process we record the time when branching starts, the time when the first solution is found, the time when a solution within 5\% of the optimal is found, and the time when the optimal solution is found. This allows us to compute the percentage of time spent on each phase, and the percentage of time spent branching versus before branching (i.e., pre-processing the instance and processing the root node). We average these numbers over the 300 samples to obtain a view of the typical behaviour of the solver on each benchmark. We further divide phase 2 (improvement) into two sub-phases: (2a) from the first feasible solution to the first feasible solution with objective value within 5\% of the optimal, and (2b) which encompasses the rest of phase 2. The results are shown in Figure \ref{fig:phase_analysis}. We observe the following. For all benchmarks, obtaining a feasible solution is trivial. For set covering instances, the optimal solution is often known by the time that branching starts. In the case of combinatorial auctions, the optimal solution is typically not known at the start of B\&B, but a good solution is. For GISP, finding optimal, or even good, solutions is not as easy, making the proving phase relatively shorter. We conclude that these benchmarks allow us to test our methodology on three very different settings that may arise in a real-life situation.

\paragraph{Data collection procedure}
For each instance, we collect information at the root node: the bipartite graph representation $X_G=(\mC, \mV, \mA)$ and the optimal root LP value $z^{LP}$. We then proceed to solve the instance. For the first 100 processed nodes and as long as no incumbent exists, no samples are collected. This allows us to initialize statistics as the trend of open nodes $\tau(t)$, and to ignore instances that are solved within 100 nodes which are therefore too easy. After 100 nodes have been processed and an incumbent exists, we collect samples with a probability of $0.02$. At sampling time, we record the value of the dynamic features (see Section \ref{sec:methodology_phase}), as well as the incumbent value $\bar{z}(t)$. Once the instance is solved, the collected samples are completed by appending the root node information $(X_G, z^{LP})$ as well as the optimal objective value $z^*$, which will be used as a target.

\paragraph{Optimal objective value prediction (Q1)}
We test the prediction accuracy of our GNN model on the four benchmarks. We train a model for each of the targets described in Section \ref{sec:methodology}. We measure the error as
\begin{equation}
\label{eq:error}
    e = 100 \times \frac{1}{N} \sum_{i=1}^N \frac{|z^*_i - \tilde{z}^*_i|}{|z^*_i|}
\end{equation}
where $N$ is the number of samples, $z^*_i$ is the optimal objective value of sample $i$ and $\tilde{z}^*_i$ is the predicted optimal objective value of sample $i$. Notice that, independently of the learning target, we measure the error in the space of the original prediction we want to make.

\paragraph{Prediction of phase transition (Q2)}
We make a prediction on whether we have transitioned to phase 3 (optimal solution has been found). We compare the performance of four predictors. The first two predictors are the ones proposed by \citet{Berthold2018}, namely $C^{\text{est}}$ (best-estimate, see Eq. \ref{eq:estimate_classifier}) and $C^{\text{rank-1}}$ (rank-1, see Eq. \ref{eq:rank1_classifier}). The third predictor $C^{GNN}_{\epsilon}$ is based on the GNN regression model, as described in Eq. \ref{eq:gnn_classifier}. We report the performance of this classifier with $\epsilon = 0$ and with a tuned value $\epsilon^*$ which was obtained by optimizing the accuracy with a small grid search over the range $[-0.02, 0.02]$ on the validation set. The fourth predictor $C^D$ is based on the dynamic features, as described in Section \ref{sec:methodology_phase}.

\subsection{Results}

Tables \ref{tab:gnn} and \ref{tab:gnn_mixed} show the results of the optimal objective value prediction task. The GNN models tested in Table \ref{tab:gnn} were trained and tested on instances of the same type. On the contrary, the results of Table \ref{tab:gnn_mixed} correspond to one unique model that was trained in the mixed dataset, and then tested on different benchmarks. First, we observe that using targets that include LP information ($\Theta_2$ and $\Theta_3$) is beneficial to performance, as opposed to directly trying to predict the optimal objective value ($\Theta_1$). There is no clear winner among targets $\Theta_2$ and $\Theta_3$. Second, we observe that the generalist model, the one trained on the mixed dataset, performs comparably to the specialized models, even outperforming them in some cases. \\

\begin{figure}
     \centering
     \begin{subfigure}[b]{0.45\textwidth}
         \centering
         \includegraphics[width=\textwidth]{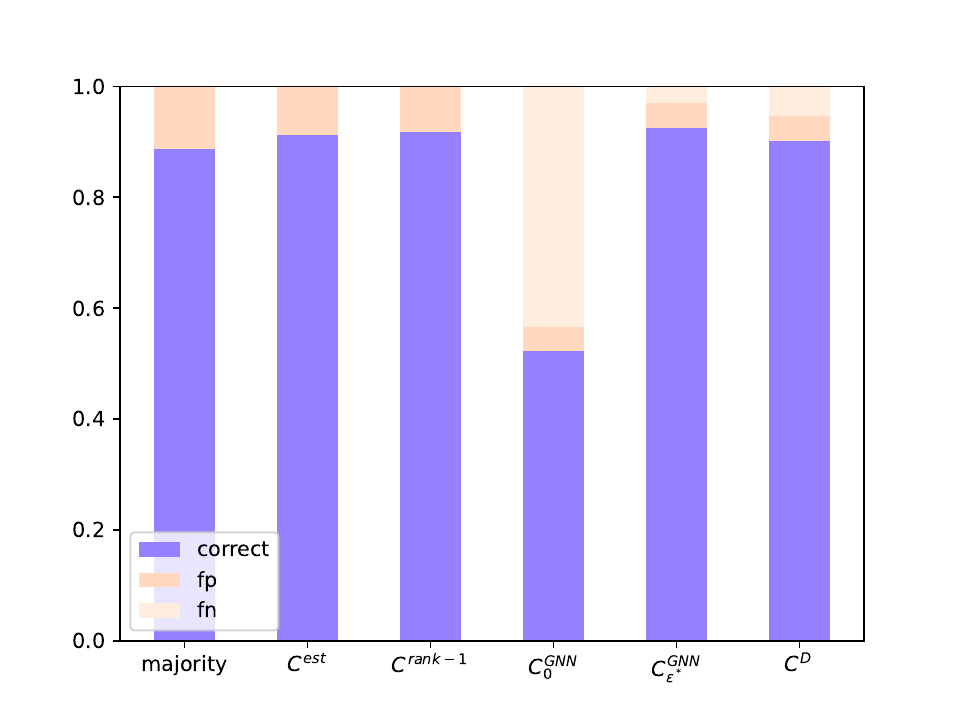}
         \caption{Set covering}
     \end{subfigure}
     \begin{subfigure}[b]{0.45\textwidth}
         \centering
         \includegraphics[width=\textwidth]{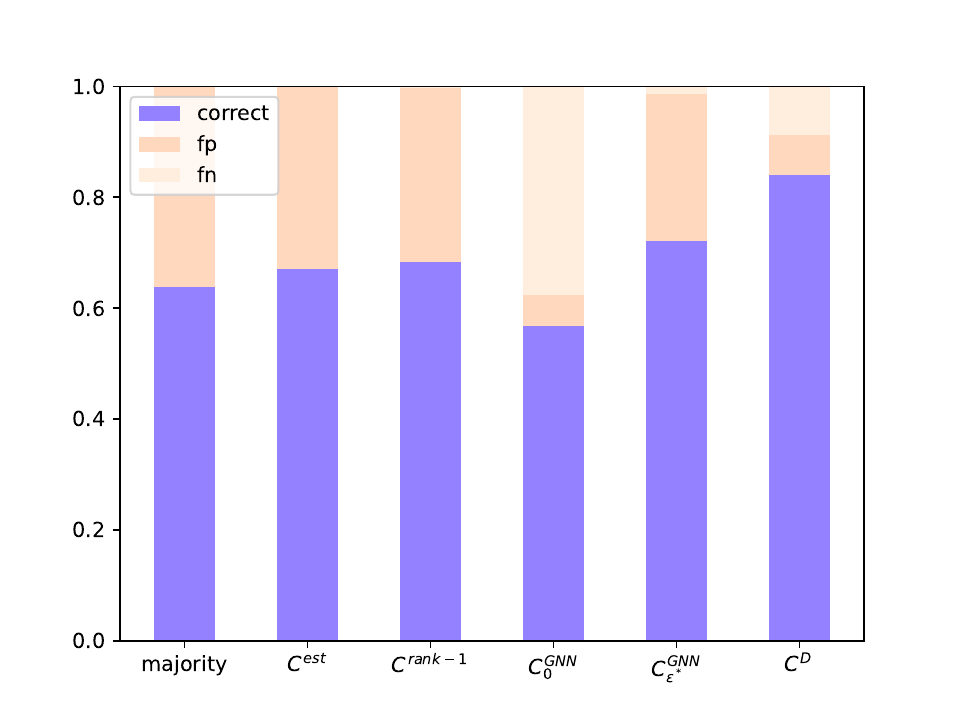}
         \caption{Combinatorial auctions}
     \end{subfigure} \\
     \begin{subfigure}[b]{0.45\textwidth}
         \centering
         \includegraphics[width=\textwidth]{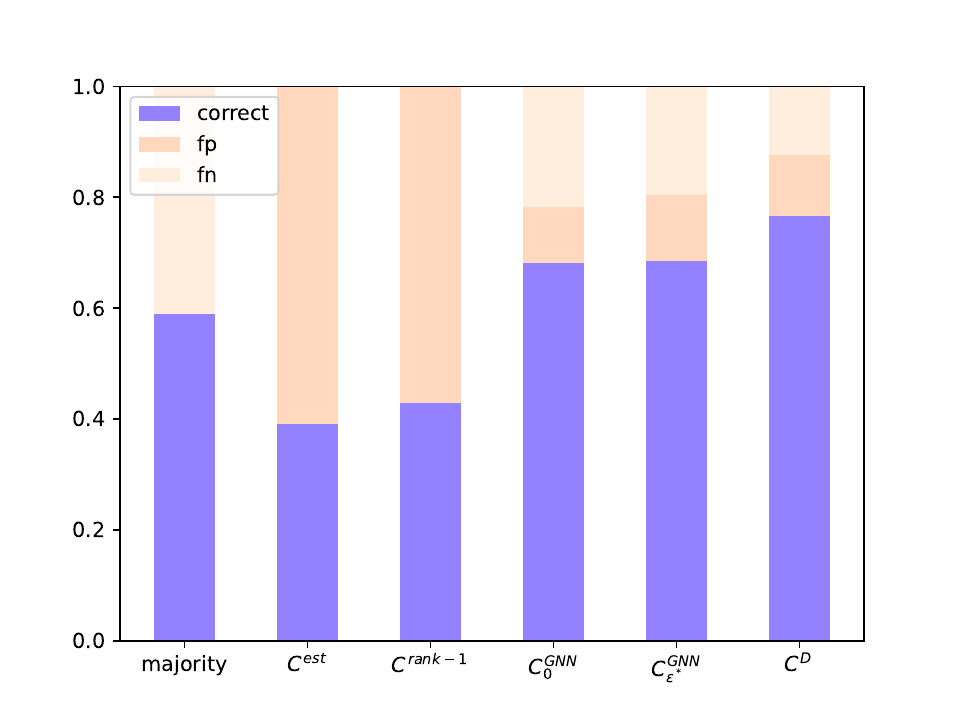}
         \caption{GISP}
     \end{subfigure}
     \begin{subfigure}[b]{0.45\textwidth}
         \centering
         \includegraphics[width=\textwidth]{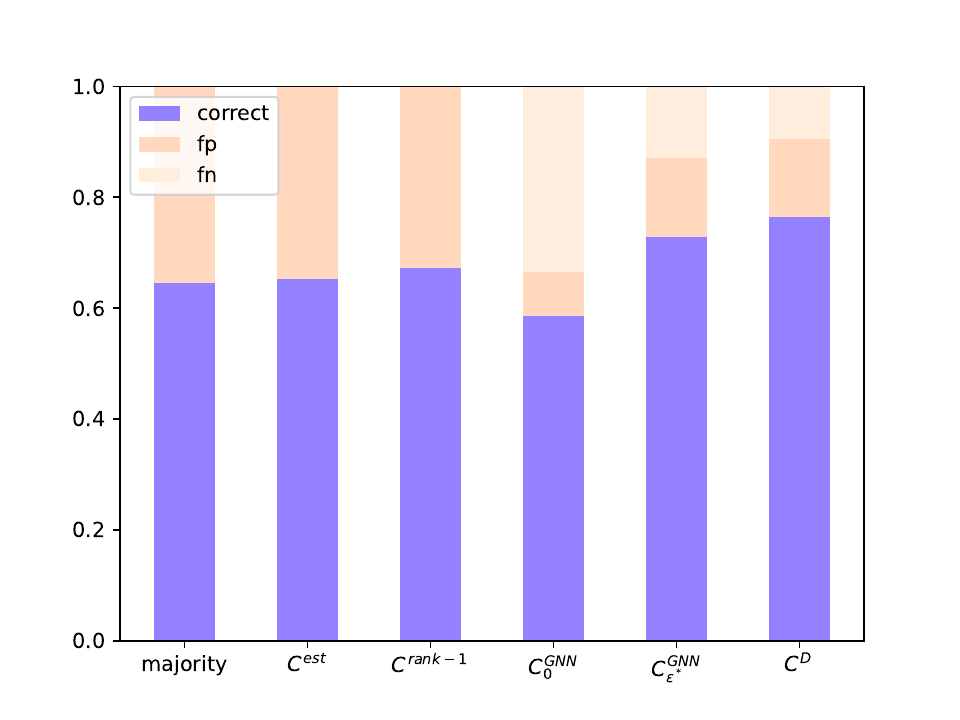}
         \caption{Mixed}
     \end{subfigure}
\caption{Prediction accuracy of the different classifier models. We show the fraction of correctly classified samples (correct, in purple), the fraction of false positives (fp, dark yellow) and the fraction of false negatives (fn, light yellow). }
\label{fig:results}
\end{figure}

We now select one GNN model per benchmark to be used in the next prediction task: the phase transition prediction. We select the model in the following way: we use the specialized model that achieves the best result on the validation set. Figure \ref{fig:results} (a-c) shows the results for all classifiers on the pure benchmarks (see Table \ref{tab:results} for the same results in table form). Further, we include a column that shows the classification accuracy of a dummy model that always predicts the majority class. We observe that the classifiers of \citet{Berthold2018} (best-estimate and rank-1) tend to predict the phase transition too early. This is, they mostly output a positive prediction, which means they believe the incumbent to be optimal. This results in the misclassified samples being almost exclusively false positives. On the contrary, the GNN model $C^{GNN}_{0}$ tends to be too pessimistic, which can be fixed with the right tuning of the $\epsilon$ parameter. For all benchmarks, $C^{GNN}_{\epsilon^*}$ performs better than the classifiers of \citet{Berthold2018}. At the same time, the inclusion of the dynamic features ($C^D$) further improves the performance, except for set covering where $C^{GNN}_{\epsilon^*}$ and $C^D$ are close to a tie. \\

It is important to notice that, depending on the application, false positives and false negatives could have very different consequences. As an example, if the phase transition prediction is used to change the behaviour of the primal heuristics (e.g. switch them off once the optimal is found) a false positive could excessively delay finding the optimal solution and therefore has a much bigger potential of harming performance than a false negative. The parameter $\epsilon$ provides an easy way to navigate this tradeoff, where one could sacrifice some accuracy to keep the rate of false positives to a minimum.\\

Figure \ref{fig:results}d shows the same experiment but on a mixed dataset. This is, the models were trained \emph{and} tested on a benchmark comprised of instances of all three types (in equal proportion). We observe a similar behaviour compared to the specialized benchmarks. The GNN model $C^{GNN}_{0}$ tends to be too pessimistic, while $C^{GNN}_{\epsilon^*}$ achieves better accuracy and better false positive rate than the classifiers of \citet{Berthold2018}. Using dynamic features further improves the accuracy of the model.\\

\begin{table*}[tb]
\centering
\caption{Average relative error (as defined in Eq. \ref{eq:error}) of the GNN model. One model was trained per benchmark. The train and test instances in each case are of the same type. }
\label{tab:gnn}
\begin{tabular}{rcccccc}
Instances & & $\Theta_1$ & & $\Theta_2$ & & $\Theta_3$ \\ \hline \hline
Set covering & & $1.48$ & & $0.80$ & & $\mathbf{0.54}$ \\
Combinatorial auctions  & & $3.20$ & & $\mathbf{0.55}$ & & $0.62$ \\
GISP & & $3.32$ & & $\mathbf{2.35}$ & & $2.39$ \\
\bottomrule[1.5pt]
\end{tabular}
\end{table*}

\begin{table*}
\centering
\caption{Average relative error (as defined in Eq. \ref{eq:error}) of the GNN \emph{mixed} model. Only one model was trained on a dataset comprised of intances of all types. The test sets are comprised of instances of one type only, except for the mixed test set (last row). }
\label{tab:gnn_mixed}
\begin{tabular}{rcccccc}
Instances & & $\Theta_1$ & & $\Theta_2$ & & $\Theta_3$ \\ \hline \hline
Set covering & & 1.35 & & \textbf{0.73} & & 0.82 \\
Combinatorial auctions & & 3.15 & & 1.17 & & \textbf{0.53} \\
GISP & & 3.17 & & \textbf{2.32} & & 2.43 \\
Mixed test set & & 1.70 & & 0.97 & & \textbf{0.75} \\
\bottomrule[1.5pt]
\end{tabular}
\end{table*}

\begin{figure}
    \centering
    \includegraphics[width=\textwidth]{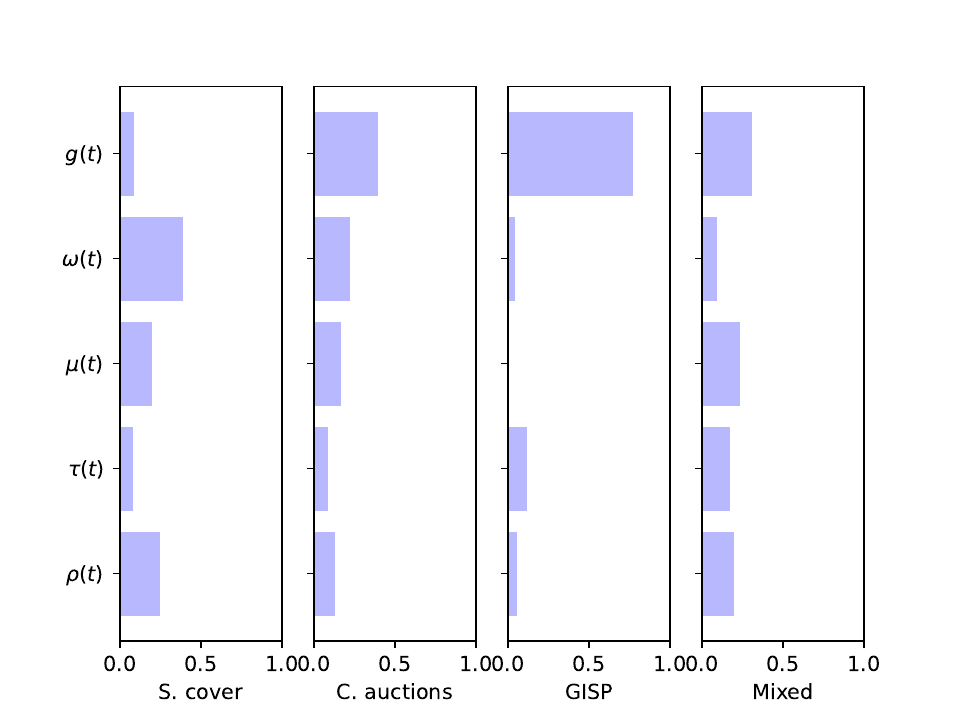}
    \caption{Feature importance of the dynamic models trained to predict phase transition for each of the benchmarks.}
    \label{fig:feature_importance}
\end{figure}

Finally, we analyze the importance of the dynamic features assigned by the $C^D$ classifier (Figure \ref{fig:feature_importance}). We see that the four learned models are in fact very different, with the GISP model mostly making decisions based on the gap and the other three considering all features more uniformly. This speaks in favour of learning on sets of instances of the same type.

\begin{table*}
\centering
\caption{Prediction accuracy of the different classifier models. We show the fraction of correctly classified samples, the fraction of false positives and the fraction of false negatives.}
\label{tab:results}
\begin{tabular}{lcccccc}
 & & Correct & & False positives & & False negatives \\ \hline \hline
 Majority & & 0.89 & & 0.11 & & 0.00 \\
$C^{\text{est}}$  & & 0.91 & & 0.09 & & 0.00 \\
$C^{\text{rank-1}}$  & & 0.92 & & 0.08 & & 0.00 \\
$C^{\text{GNN}}_0$ && 0.52 && 0.05 && 0.43 \\
$C^{\text{GNN}}_{\epsilon^*}$ && 0.93 && 0.04 && 0.03 \\
$C^{D}$ && 0.90 && 0.05 && 0.05 \\ \bottomrule[1.5pt]
\multicolumn{7}{c}{Set covering} \\ \\ 
 & & Correct & & False positives & & False negatives \\ \hline \hline
 Majority & & 0.64 & & 0.36 & & 0.00 \\
$C^{\text{est}}$  & & 0.67 & & 0.33 & & 0.00 \\
$C^{\text{rank-1}}$  & & 0.68 & & 0.32 & & 0.00 \\
$C^{\text{GNN}}_0$ && 0.57 && 0.06 && 0.37 \\
$C^{\text{GNN}}_{\epsilon^*}$ && 0.72 && 0.27 && 0.01 \\
$C^{D}$ && 0.84 && 0.07 && 0.09 \\ \bottomrule[1.5pt]
\multicolumn{7}{c}{Combinatorial auctions} \\ \\ 
 & & Correct & & False positives & & False negatives \\ \hline \hline
Majority & & 0.59 & & 0.00 & & 0.41 \\
$C^{\text{est}}$  & & 0.39 & & 0.61 & & 0.00 \\
$C^{\text{rank-1}}$  & & 0.43 & & 0.57 & & 0.00 \\
$C^{\text{GNN}}_0$ && 0.68 && 0.10 && 0.22 \\
$C^{\text{GNN}}_{\epsilon^*}$ && 0.69 && 0.12 && 0.19 \\
$C^{D}$ && 0.77 && 0.11 && 0.12 \\ \bottomrule[1.5pt]
\multicolumn{7}{c}{GISP} \\ \\
 & & Correct & & False positives & & False negatives \\ \hline \hline
Majority & & 0.64 & & 0.36 & & 0.00 \\
$C^{\text{est}}$  & & 0.65 & & 0.35 & & 0.00 \\
$C^{\text{rank-1}}$  & & 0.67 & & 0.33 & & 0.00 \\
$C^{\text{GNN}}_0$ && 0.59 && 0.08 && 0.34 \\
$C^{\text{GNN}}_{\epsilon^*}$ && 0.73 && 0.14 && 0.13 \\
$C^{D}$ && 0.77 && 0.14 && 0.09 \\ \bottomrule[1.5pt]
\multicolumn{7}{c}{Mixed} \\
\end{tabular}
\end{table*}

\section{Conclusions}
\label{sec:conclusions}

In this paper, we presented our methodology for predicting the optimal objective value of MILPs. Compared to the literature on predicting optimal solutions, our learning task is easier, yet still offers a variety of possibilities for its application within MILP solvers. Our methods can be used to both predict the optimal objective value and to classify a feasible solution into optimal or sub-optimal. Our computational study shows that our proposed approach outperforms the existing approaches in the literature. Further, they provide more flexibility to tune the model into the desired behaviour. 
We show that there are benefits to learning a model that specializes to an instance type, yet our model is still able to generalize well and have superior performance to other methods on mixed instance sets. 

These results open the door for many possible applications. In general terms, this prediction can be used to adapt the behaviour of the different solver components and rules depending on the solving phase. These applications, however, require further study and will be the subject of future work.

\clearpage
\bibliographystyle{abbrvnat}
\bibliography{references}

\end{document}